\documentclass[12pt]{article}
\usepackage{latexsym,amssymb,amsmath,amsfonts,amsthm,graphicx,url}
\usepackage[notref,notcite]{showkeys}
\theoremstyle{definition}

\def\beq{ \begin{equation} }
\def\eeq{ \end{equation} }
\def\mn{\medskip\noindent}

\def\square{\vcenter{\vbox{\hrule height .4pt
  \hbox{\vrule width .4pt height 5pt \kern 5pt
        \vrule width .4pt} \hrule height .4pt}}}

\def\ZZ{\mathbb{Z}}

\def\clearp{\clearpage}

\begin{document}

\title{Triple birthday matches in the Senate:\\
Lies, damned lies, and chatGPT
\thanks{The second line of the title is a famous reamrk of Samuel Clemens about statisitcs. }}
\author{Rick Durrett, Duke University \footnote{Comments, corrections, and insults should be sent to rtd@math.duke.edu} }

\date{\today}						

\maketitle

\begin{abstract}
Our question is ``What is the probability that at least three members of the senate share the same birthday?'' Before the pandemic, I asked this question in several popular math talks I gave at universities across the country. Inspired by ChatGPT's abysmal failure to answer the question, I have recently come back to this problem and now have a more satisfactory answer, thanks in no small part to what I learned form a page of Wolfram's Math World, which I located by a Google search. 
\end{abstract}

\section{Introduction}

I begin with the answer I gave in a {\it Math for Everyone} talk at Notre Dame in April 2018

\begin{itemize}

\item
The number of senators born on a given day is Binomial(100,1/365) 
and hence approximately Poisson with mean 100/365.

\item
Probability three senators born on a given day
$$
e^{-100/365} (100/365)^3/3! = 0.002606
$$
The expected number of triple birthdays in 365 days $= 0.9512$. If we assume that the events $B^3_j =$ exactly three senators were born on day $j$ are independent then the number would be Poisson with mean $0.9512$ and the probability that this occurs for at least one day is $1-\exp(-0.9512) =0.6137$

\item
Probability two senators born on a given day  
$$
e^{-100/365} (100/365)^2/2! = 0.028536
$$
expected number among 365 days = 10.415. Again if we assume that the events $B^2_j =$ exactly two senators were born on day $j$ are independent then the number would be Poisson with mean 10.415.

\item
 Answers for the 2018 Senate: One triple birthday, 
 May 3: Jim Risch (Idaho), David Vitter (Louisiana), Ron Wyder (Oregon).
Ten double birthdays: Jan 7, March 31, June 22, August 24, Sept 29,
October 20 and 24, November 17, December 7 and 10

\end{itemize}

\medskip
Curioius at what chatGPT would say I typed in:

\mn
What is the probability that at least three members of the senate share the same birthday?

\mn
The answer it gave was: 

\begin{quote}
The probability that at least three members of a senate of size n share the same birthday is approximately 
$$
1 - (364/365)^{\binom{n}{3}} - (364/365)^n + (364/365)^{\binom{n}{2}}.
$$
This is because there are 365 possible birthdays, and the chance of any one senator not sharing a birthday with any other is 364/365. The chance of no two senators sharing a birthday is $(364/365)^{\binom{n}{2}}$, and the chance of no three senators sharing a birthday is $(364/365)^{\binom{n}{3}}$. So, to get the chance of at least three senators sharing a birthday, we subtract the chance of no two sharing a birthday, no three sharing a birthday and add the chance of all of them sharing a birthday. For example, if there are 100 members in a senate, the probability that at least three members share the same birthday is approximately 0.527 or 52.7\%.
\end{quote}

\medskip
To start to understand this answer we consider the claim: ``The chance of no two senators sharing a birthday is $(364/365)^{\binom{n}{2}}$.'' This would be true if the events $A_{i,j}$ that the $i$th and $j$th senators have the same birthday were independent but they are not. If $\{i,j\} \neq \{k,\ell\}$ are distinct integers then 
$$
P(A_{i,j} \cap A_{k,\ell}) = \frac{1}{365^2} = P(A_{i,j}) P(A_{k,\ell})
$$
This shows that the events are {\bf pairwise independent}, but for true independence we would need
$$
P(A_{1,2} \cap A_{2,3} \cap A_{3,1}) = P(A_{1,2}) P(A_{2,3}) P(A_{3,1})
$$
Unfortunately  the left-hand side is $1/(365)^2$ while the right is $1/(365)^3$. A second problem with this part of chatGPT's answer is that it is a well-known fact that the chance of no two senators sharing a birthday is
\beq
\frac{ 365 \cdot 364 \cdots 266}{(365)^{100}} = 3.072 \times 10^{-7}
\label{ans1}
\eeq
versus $(354/365)^{\binom{100}{2}} = 1.265 \times 10^{-6}$

The second claim ``the chance of no three senators sharing a birthday is $(364/365)^{\binom{n}{3}}$'' is more mysterious.
The try to explain this we turn to Suraj Regmi's blog post from January 26, 2019, which, as this was being written, was top answer in response to the Google search for ``birthday triples math formulas.'' Regmi reasons: There are $C(n,3)$ triplets and for all of the triplets to not have the same birth date probability becomes
\beq
( 1 - 1/(365)^2)^{C(n,3)}
\label{ans2}
\eeq
as the triples are INDEPENDENT EVENTS. 

While the situation with double birthdays is subtle, involving the distinction between pairwise and fully independent events, the situation for triples is not. If $A_{i,j,k}$ is the probability that senators $i,j,k$ share the same birthday then
$$
P( A_{1,2,3} \cap A_{2,3,4} ) = \frac{1}{(365)^3} < 
 \frac{1}{(365)^4} = P( A_{1,2,3}) \cdot P( A_{2,3,4} ) 
$$
When $n=100$ the formula in \eqref{ans2} evaluates to 0.29708, which gives an answer of \beq
1-0.29708 = 0.70292
\label{ans3}
\eeq

To be fair to the bot and the blogger, I should confess that I sinned at Notre Dame: my events are not independent
\begin{align*}
P(B^2_1) = 0.037155 &\qquad P(B^2_2|B^2_1) = 0.035676 \\
P(B^3_1) = 0.003325 &\qquad P(B^3_2|B^3_1) = 0.003032
\end{align*}
However, I did not claim them to be. My goal was to give a simple approximation for the probability that there is a day on which exactly three senators were born and it turns out to be quite accurate. When you put 100 in the Wolfram's birthday problem calculator \cite{WolfA} you get the answers for birthdays in senate.

\begin{center}
\begin{tabular}{ll}
at least 2 the same & $1 - 3.072 \times 10^{-7}$ \\
at least 3 the same & 0.6459
\end{tabular}
\end{center}

\mn
The first result is 1 minus \eqref{ans1}. The second answer shows that 
\eqref{ans2} and \eqref{ans3} are wrong. 

Note that here 0.6459 is the probability for at least 3 senators sharing a birthday, compared to the earlier estimate of 0.6137 for exactly 3. At the end of Section 2.1 the reader we will see that our answer using the Poisson approximation is very close to the answer of 0.6140 for exactly 3 senators sharing a birthday computed by using the first six terms of the inclusion exclusion formula.  Our approximation for the expected number of double birthdays 10.415 compares very well  with the true expected value of 10.3645, which should not be surprising since expected values are not sensitive to dependence. That fact is fortunate for us, since our argument implies that the number of double birthdays has a Poisson distribution, but as results below will show (see Figure 1) this statement is not very accurate. 

\clearpage

\section{Calculations}

\subsection{The probability of a triple birthday}

Let $T$ be the number of triple birthdays, i.e.,  the days of the year that are the birthdays of exactly three senators. If we use $i$ as shorthand for $(i_1,i_2,i_3)$ with $1\le i_1 < i_2 < i_3 \le 100$ and let $A_i$ be the event that senators $i_1,i_2,i_3$ have the same birthday and no other senator has this birthday then
\begin{align*}
q_1= \sum_i P(A_i) & = \binom{100}{3} \left( \frac{1}{365}\right)^2 
 \left( \frac{364}{365}\right)^{97} \cdot \frac{365}{365}  \\
 & = \binom{100}{3} \frac{365 (364)^{97}}{(365)^{100}}= 0.9301
\end{align*}
which gives $ET$ and an upper bound on the probability of $\cup_i A_i$. To get a lower bound using the Bonferroni inequalities, we need to subtract
\begin{align*}
q_2 = \sum_{i<j} P(A_i \cap A_j) & =
\frac{1}{2}  \cdot \binom{100}{3} \left( \frac{1}{365}\right)^2 
 \binom{97}{3} \frac{364}{365} \left( \frac{1}{365}\right)^2 
 \left( \frac{363}{365}\right)^{94} \\
& =
\frac{1}{2} \cdot  \binom{100}{3}  \binom{97}{3} 
\frac{365 \cdot 364 \cdot (363)^{94}} {(365)^{100}}
= 0.3996
\end{align*}
 (here $<$ is lexicographic or dictionary order on triples $(i_1,i_2,i_3)$).

To get a second upper bound we need to add
\begin{align*}
q_3 & = \sum_{i<j<k} P(A_i \cap A_j \cap A_k) \\
& =\frac{1}{3!} \cdot  \binom{100}{3}  \binom{97}{3} \binom{94}{3} 
\left( \frac{1}{365}\right)^6 \frac{364}{365}\cdot   \frac{363}{365}
\left( \frac{362}{365}\right)^{91} \\
& =\frac{1}{3!} \cdot  \binom{100}{3}  \binom{97}{3} \binom{94}{3} 
\frac{P_{365,3} \cdot (362)^{100-9}} {(365)^{100}}
= 0.1054
\end{align*}
To get a second lower bound, we need to subtract  
$$
q_4  = \sum_{i<j<k<\ell} P(A_i \cap A_j \cap A_k \cap A_\ell)
$$
To do this we use the general formula
$$
q_k =\frac{1}{k!} \prod_{k=0}^{k-1} \binom{100-3j}{3} \cdot
\frac{P_{365,k} \cdot (365-k)^{100-3k}} {(365)^{100}}
$$
which gives $q_4=0.019153181$, $q_5 =2.548039 \times 10^{-3}$, and
$q_6 = 2.57641 \times 10^{-4}$

To computer the answer 

\begin{center}
\begin{tabular}{cr}
upper bound $u_1$ & $q_1=0.931045$ \\
lower bound $v_1$ & $q_1-q_2 = 0.530545$ \\
upper bound $u_2$ & $v_1+q_3 = 0.635962$ \\
lower bound $v_2$ & $u_2-q_4 = 0.616809$ \\
upper bound $u_3$ & $v_2+q_5 = 0.614261$ \\
lower bound $v_3$ & $u_3-q_6 = 0.614004$
\end{tabular}
\end{center}

\mn
Later we will need to do this calculation for $n$ people and a calendar with $d$ days. In this case the $k$th term in inclusion-exclusion is
$$
q_k(n,d) =\frac{1}{k!} \prod_{k=0}^{k-1} \binom{n-3j}{3} \cdot
\frac{P_{d,k} \cdot (d-k)^{100-3k}} {(365)^{100}}
$$

\subsection{The number of double birthdays}

Write $k$ as shorthand for $(k_1,k_2)$ with $1\le k_1 < k_2\le 100$. Let $C_k$ be the event that senators $k_1$ and $k_2$ have the same birthday which is not shared by any of the other members of the Senate, and let $D$ be the number of double birthdays. Using $\binom{100}{2}=4950$
\beq
ED = \sum_k P(C_k) = 
\binom{100}{2} \left( \frac{1}{365}\right)  \left( \frac{364}{365}\right)^{98}
= 10.3645
\label{Edouble}
\eeq
in contrast to the approximate value of 10.415. The error comes from the Poisson approximation of the binomial. 
$$
365 \cdot P( \hbox{binomial}(100,1/365)=2 )
365 \cdot 0.028396 = 10.3645
$$

If $D$ had a Poisson($\lambda$) distribution then $ED(D-1)=\lambda^2$. Writing $j<k$ for the lexicographic order on $\ZZ^2$ and noting that birthday coincidences are pairwise independent
\begin{align*}
ED(D-1) &= \sum_{j<k} P(C_j \cap C_k) \\
& = \binom{100}{2} \left( \frac{1}{365}\right)  
\binom{98}{2}  \left( \frac{364}{365} \cdot \frac{1}{365}\right) 
\left( \frac{363}{365}\right)^{96} \\
& = ED \cdot \binom{98}{2}  \left( \frac{1}{364}\right)  \left( \frac{363}{364}\right)^{96} = ED \cdot 10.027 < (ED)^2
\end{align*}
where in the last step we have multiplied and divided by $(364/365)^{98}$

Hocking and Schweterman \cite{HS86} have derived a formula for the probability $p_k$ of $k$ double birthdays and no triple (or higher) coincidences. As we have already noted in \eqref{ans1}
$$
p_0 = \frac{P_{365,100}}{(365)^{100}} =  3.072 \times 10^{-7}
$$
Arguing as in our treatment of triple birthdays
\begin{align*}
p_1 & = \binom{100}{2} \frac{1}{365} \cdot \frac{P_{364,98}}{(365)^{98}}
\cdot\frac{365}{365}
= \binom{100}{2} \cdot \frac{P_{365,99}}{(365)^{100}}\\
p_2 & =  \frac{1}{2} \binom{100}{2} \frac{1}{365} \cdot \binom{98}{2} \frac{364}{365} \cdot  \frac{1}{365} \cdot  \frac{P_{363,96}}{(365)^{96}}\\
&=  \frac{1}{2} \binom{100}{2}\binom{98}{2}  
\cdot  \frac{P_{365,98}}{(365)^{100}}\\
p_3 & =  \frac{1}{3!} \binom{100}{2} \frac{1}{365} 
\cdot \binom{98}{2} \frac{364}{365} \cdot  \frac{1}{365} 
\cdot \binom{96}{2} \frac{363}{365} \cdot  \frac{1}{365} 
\cdot  \frac{P_{362,94}}{(365)^{94}} \\
& =  \frac{1}{3!} \binom{100}{2} \binom{98}{2} 
\binom{96}{2} \cdot  \frac{P_{365,97}}{(365)^{100}} 
\end{align*}
Referring to (1) in \cite{HS86} and doing some algebra, we see that the formula for general $k$ and $n$ is
$$
p_k  = \frac{1}{k!} \prod_{j=0}^{k-1} \binom{n-2j}{2} \cdot 
\frac{P_{365,n-k}}{(365)^n} 
$$
Using this formula we see that
$$
P_{k} = P_{k-1} \cdot \frac{1}{k} \cdot  \binom{n-2(k-1)}{2} \cdot \frac{1}{365-n+k}
$$
If $p_k$ was $\alpha$ times the Poisson($\lambda$) distribution (recall that $\sum_k p_k = P(T=0)$) then we would have $p_k = \lambda p_{k-1}$ so $p_k/P(T=0)$ is not a Poisson distribution. To compare with the Poisson (see Figure 1) we note that
\beq
\sum_k k P_k = 3.87454 \qquad
\sum_k p_k =P(T=0) = 0.354135
\label{notrip}
\eeq
which agrees with the Wolfram Alpha result $P(T=0)=1-0.6549$. So we have $E(D|T=0) = 10.941$.

\begin{figure}[h]
  \centering
  \includegraphics[width=4in,keepaspectratio]{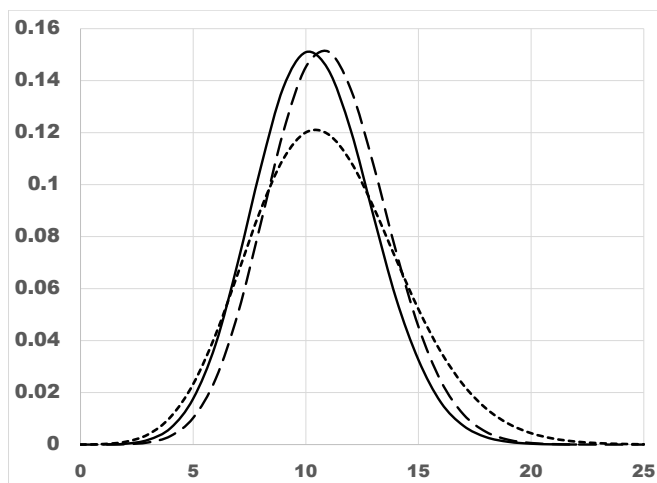}
\caption{Graph of $p_k =$ the probability of $k$ double birthdays conditioned on no triple birthday (line with longer dashes), compared to Poisson with mean 10.941 (shorter dashes). Solid line gives results from simulation of 1 million instances of the unconditioned distribution. The mean is 10.36 in agreement with \eqref{Edouble}. }
\end{figure}

\subsection{The number of triple birthdays}

McKinney \cite{McK}was perhaps the first person to try to determine the probability that in $n$ people selected at random $r$ will have the same birthday. To do this he let $X_i$, $1 \le i \le n$ be i.i.d.~uniform on $\{1, 2, \ldots M\}$. Let $n_i$ be the number of values that appear $i$ times in the sample. His main result is
\beq
P_n(n_1, n_2, \ldots n_{r-1}) = \frac{n!}{\prod_{i=1}^{r-1} n_i! (j!)^{n_j}}
\frac{ P_{M,\sum_{i=1}^{r-1} n_i}} {M^n}
\label{McKf}
\eeq

\begin{proof}
The second factor in \eqref{McKf} is the probability than $n$ independent uniforms will have $n_1$ nonrepeated values, $n_2$ pairs, $n_3$ triples in a specified order. The first factor representa the number of distinguishable ways that this particular oder can be permuted.
\end{proof}

\mn
If we let $G_r^c =$ no value is repeated $r$ or more times in the sample then
$$
P(G_r^c) = \sum\{  P_n((n_1, n_2, \ldots n_{r-1}) : \sum_i in_i =n \}
$$
He used this to compute the following values of $P(G+r)$. (His $E=G_r^c$.)

\begin{center}
\begin{tabular}{ccccc}
$r=2$ & $n=22$ & 0.4758 & $n=23$ &  0.5074 \\
$r=3$ & $n=87$ & 0.4998 & $n=88$ &  0.5114 \\
$r=4$ & $n=186$ & 0.4758 & $n=187$&  0.5033
\end{tabular}
\end{center}

To compute the distribution of the number of triple birthdays, we write recursions that are inspired by those given in Wolfram's Math World \cite{WolfMW} for $Q_i(n,d) =$ the probability that in a group of size $n$ with $d$ possible birthdays, a birthday is shared by extactly $i$ (and no more) people. Here we let $\tau_k(n,d)$ be the probability that there are exactly $k$ triple birthdays in a group of size $n$ with $d$ possible birthdays.

We have computed $\tau_0(100,365) = 0.386$.
$$
\tau_1(100,365) = \binom{100}{3} \frac{1}{365^2} 
\left( \frac{364}{365} \right)^{97}\times  \tau_0(97,364)
$$
The term $(364/365)^{97}$ is the probability that the 97 remaining people do not have birthdays that match the triple birthday. If we condition on this event then their birthday are uniform over the remaining 364 possibilities. Similarly
\begin{align*}
\tau_2(100,365) &= \frac{1}{2!} \binom{100}{3} \frac{1}{365^2} 
\binom{97}{3} \frac{364}{365} \frac{1}{365^2}
\left( \frac{363}{365} \right)^{94}\times  \tau_0(94,363)\\
&= \frac{1}{2!} \binom{100}{3} \binom{97}{3} \cdot
\frac{ P_{365,2} \cdot (363)^{94}}{(365)^{100} }\times \tau_0(94,363)
\end{align*} 
Following the pattern we can see
\begin{align*}
\tau_3(100,365) &= \frac{1}{3!} \binom{100}{3} \binom{97}{3} \binom{94}{3}\cdot
\frac{ P_{365,3} \cdot (362)^{91}}{(365)^{100} }\times \tau_0(91,362) \\
\tau_4(100,365) &= \frac{1}{4!} \prod_{j=0}^3 \binom{100-3j}{3} \cdot
\frac{ P_{365,4} \cdot (361)^{88}}{(365)^{100} }\times\tau_0(88,361)
\end{align*}

These probabilities are easier to compute than one might expect. The quantities to the left of the $\times$ signs are the $q_k$ computed in Section 2.1. So it remains to compute $\tau_0(100-3k,365-k)$ using the Bonferroni inequalities (and stopping with the fourth bound). The next table gives the results and compare with values computed from 1 million simulations. 

\begin{center}
\begin{tabular}{ccclc}
$k$ & $q_k$ & $1-\tau_0(100-3k,365-k)$ & calculation & simulation\\
0 &&& 0.386 & 0.380921\\
1 & 0.93014 & 0.58796 & 0.38325 & 0.381977\\
2 & 0.39960 & 0.55777 & 0.17672 & 0.176321\\
3 & 0.10542 & 0.52719 & 0.049843 & 0.049634\\
4 & 0.019153 & 0.49585 & 0.009656 & 0.009604 \\
5 & $2.548 \times 10^{-4}$ & 0.46415 & 0.001365 & 0.001375
\end{tabular}
\end{center}

\section*{Acknowledgement}
This work was partially supported by NSF grant DMS 2153429 from the probability program. The views expressed here are those of the author and do not necessarily reflect the view of the Natioanl Science Foundation. Computations were performed by my student Hwai-Ray Tung, who will graduate from Duke in May 2023 and go to a postdoctoral position in the Utah math department on July 1. On that date his adviser will become a James B. Duke Emeritus Professor of Mathematics.

\clearp

\end{document}